\documentclass[10pt,reqno]{article}
\usepackage{amssymb,amsmath,amsthm}
\usepackage[pdftex]{graphicx}
\usepackage{hyperref}
\usepackage{mathrsfs}
\usepackage{MnSymbol}
\hypersetup{pdfpagemode=FullScreen,colorlinks=true}

\title{H\"older maps from Euclidean spaces to Carnot groups}
\author{Z. Balogh, A. Kozhevnikov and P. Pansu\footnote{Z.B. and P. P.~are supported by MAnET Marie Curie Initial Training Network. P. P.~is supported by Agence Nationale de la Recherche grants ANR-2010-BLAN-116-01 GGAA and ANR-15-CE40-0018 SRGI, by ESPRC under Isaac Newton Institute's grant EP/K032208/1, and by Simons Foundation.
}}

\newtheorem{thm}{Theorem}
\newtheorem{prop}[thm]{Proposition}
\newtheorem{lem}[thm]{Lemma}

\newtheorem{cor}[thm]{Corollary}
\newtheorem{defi}[thm]{Definition}
\newtheorem{rem}[thm]{Remark}
\newtheorem{exa}[thm]{Example}

\newenvironment{pf}{\begin{trivlist}\item[]{\bf Proof\ }}
{\mbox{}\hfill\rule{.08in}{.08in}\end{trivlist}}

\def\R{\mathbb{R}}

\begin{document}
\maketitle

\abstract{We give alternative proofs of (unsharp) results of Gromov's on his H\"older equivalence problem: for which $\alpha$ does there exist a $C^\alpha$-homeomorphism of an open set of Euclidean space to an open set of a given Carnot group? We indicate a possible route to sharp bounds.}

\tableofcontents

\section{Introduction}

\subsection{Gromov's H\"older equivalence problem}

More than 20 years ago, M. Gromov advertised the following question. Let $G$ be a Carnot group of dimension $n$. For which $\alpha$ does there exist a $C^\alpha$-homeomorphism of an open set of $\R^n$ to an open set of $G$ ?

\medskip

In \cite{Gromovcarnot}, Gromov provides two approaches to this problem. The first consists in proving lower bounds on the Hausdorff dimensions of subsets of a given topological dimension. The second uses norms on Alexander-Spanier representatives of cohomology classes.

A third approach has been provided by Roger Z\"ust. It consists of integrating differential forms along H\"older maps. In this note, we propose a variant of Z\"ust's method that provides an alternative proof of most of Gromov's numerical bounds, \cite{Gromovcarnot}. We state the result for the larger class of Carnot \emph{manifolds}, defined in subsection \ref{mani} below.

\begin{thm}
\label{thm}
Let $M$ be an $n$-dimensional equiregular and equihomological Carnot manifold, equipped with a Carnot-Carath\'eodory metric. Let $\alpha(G)$ be the supremum of $\alpha<1$ such that there exists a $C^\alpha$-homeomorphism of an open subset of Euclidean $n$-space on an open subset of $M$. 
\begin{enumerate}
  \item If $M$ has Hausdorff dimension $Q$, then $\alpha(M)\leq \frac{n-1}{Q-1}$.
  \item If $n=4m+3$ and $G$ is the $m$-th quaternionic Heisenberg group, then $\alpha(G)\leq \frac{4m+1}{4m+4}$.
  \item If $n=2m+1$ and $M$ is a contact manifold, then $\alpha(M)\leq \frac{m+2}{m+1}$.
  \item If $M$ is a generic distribution of rank $h$, and if $h-k\geq (n-h)k$, then $\alpha(M)\leq \frac{n-k}{n+h-k}$.
\end{enumerate}
\end{thm}

We are unable to attain sharp bounds. For instance, for a step $s$ Carnot group, the group exponential $\mathfrak{g}\to G$ is a $C^{1/s}$-homeomorphism, we expect this to be optimal when $s=2$, but we are unable to prove it yet. Nevertheless, in section \ref{sharp}, we describe how a sharp bound would follow (at least in dimension 3) from a conjectural slicing theorem in a theory of currents adapted to Carnot groups.

\subsection{Strategy}

It is inspired by R. Z\"ust's work, \cite{Zust}, \cite{LDZ}, \cite{Zusttrees}, \cite{Zust2016}. Z\"ust defines the integral of a smooth differential form on $G$ along a $C^\alpha$-H\"older continuous map, in such a way that Stokes' Theorem holds. Homogeneity of differential forms commands the exponent $\alpha$. When $\alpha$ is too large, all integrals of $k$-forms of high weight vanish. This indicates that the induced morphism on  homology vanishes, which is unlikely for a homeomorphism. This contradiction concludes the proof.

We choose to modify Z\"ust's construction slightly and define integrals as limits of dyadic PL approximations.

\subsection{Plan of the paper}

Section \ref{subd} starts with an interpolation procedure for continuous maps of cubes to Carnot groups. In section \ref{integral}, sufficient conditions are given on the H\"older exponent in order that smooth differential forms on a Carnot group can be integrated along H\"older continuous maps from the cube, and when such integrals have to vanish. This yields upper bounds on $\alpha(G)$ in section \ref{bounds}. Section \ref{sharp} explains a possible route to the sharp bound $\frac{1}{2}$ for Heisenberg group $\mathbb{H}^1$.

\subsection{Acknowledgements}

Thanks to Piotr Hailasz for is interest in this work, to Gabriel Pallier for his explanations concerning slicing, and to Roger Z\"ust for inspiring discussions.

\section{Subdivisions}
\label{subd}

We define a sequence of continuous, piecewise smooth approximations of a H\"older continuous maps from a subset $Q$ of $\R^k$ to $G$. For this, we subdivide $Q$. Since we want pieces which keep the same shape, it is easier to subdivide cubes than simplices. On the other hand, it is easier to interpolate a map on simplices than on cubes. Therefore one needs to be able to pass from cubes to simplices, i.e. to triangulate cubes.

\subsection{Straight simplices in Carnot groups}
\label{straight}

Let $\Delta^k$ denote the standard $k$-simplex. Given a $k+1$-tuple of points $(p_0,\ldots,$ $p_k)\in G^k$, define the \emph{straight simplex} $\sigma_{(p_0,\ldots,p_k)}:\Delta^{k}\to G$ inductively as follows: construct the affine cone with vertex at 0 on the map $\log\circ L_{p_0^{-1}}\circ\sigma_{(p_1,\ldots,p_k)}):\Delta^{k-1}\to\mathfrak{g}$, apply $\exp$ and left translate by $p_0$. When $G$ is 2-step, left translations are affine in exponential coordinates, thus the straight simplex coincides with the affine simplex in exponential coordinates, whence the analogy with PL maps. For higher step groups, the straight simplices need not be affine in any sense any more. Nevertheless, $\sigma_{(p_0,\ldots,p_k)}$ is smooth and there exists a constant $C$ such that
\begin{eqnarray*}
\mathrm{diameter}(\{p_0,\ldots,p_k\})\leq 1 \Rightarrow \mathrm{diameter}(\sigma_{(p_0,\ldots,p_k)})\leq C\,\mathrm{diameter}(\{p_0,\ldots,p_k\}).
\end{eqnarray*}

\subsection{Triangulating the cube}

Let $Q=[0,1]^k$ be the unit $k$-cube. Index vertices of $Q$ by strings in $\{0,1\}^k$. Any injective map $s:\{0,1,\ldots,k\}\to\{0,1\}^k$ defines a linear simplex $\sigma_s$ (affine map of the standard simplex to $\R^k$) with vertices at the corners of the cube. Its orientation is given by the sign of the determinant
\begin{eqnarray*}
\mathrm{sign}(\mathrm{det}(s(1)-s(0),\ldots,s(k)-s(0))),
\end{eqnarray*}
which we denote by $(-1)^s$.

There is an obvious partial order on strings: $x\leq y$ if for all $i$, $x_i\leq y_i$. Denote by $\mathcal{I}$ the set of increasing maps $\{0,1,\ldots,k\}\to\{0,1\}^k$. The images of simplices $\sigma_s$, $s\in\mathcal{I}$ constitute a triangulation of $Q$. The singular chain
\begin{eqnarray*}
\sum_{s\in\mathcal{I}}(-1)^s \sigma_s
\end{eqnarray*}
represents the fundamental class of $Q$. 

\subsection{Interpolation}
\label{interpolation}

Let $F:\{0,1\}^k\to G$ be a map defined on the vertices of the cube. Let $F':Q\to G$ be the map such that for every $s\in\mathcal{I}$, $F'\circ \sigma_s=\sigma_{(s(0),\ldots,s(k))}$. Then $F'$ is continuous. Indeed, faces of simplices match, since vertices occurring in two different simplices appear in the same order, and hence define the same simplex in both.

Let $F:Q\to G$ be a continuous map. Consider the regular subdivision of the unit cube in $2^{jk}$ tiny cubes. Apply the above procedure to each tiny cube. The resulting map $F_j:[0,1]^k\to G$ (made of $k!2^{jk}$ pieces) is the required piecewise smooth interpolation. It is continuous. Indeed, if $(x_0,\ldots,x_\ell)$ is a sequence of vertices belonging to two tiny cubes $Q_0$ and $Q_1$ adjacent along the $i$-th face, their labels as vertices of $Q_1$ are obtained from their labels as vertices of $Q_0$ by switching the 0 at the $i$-th place to 1. If the sequence is increasing in $Q_0$, it is increasing in $Q_1$ as well, so the straight simplices used by both tiny cubes coincide. It follows that smooth differential forms can be pulled back by $F_j$ to piecewise smooth forms which satisfy Stokes Theorem.

As $j$ tends to $\infty$, $F_j$ converges uniformly to $F$.

\subsection{Passing from $F_j$ to $F_{j+1}$}

Let $\{0,1,2\}^k$ be the set of vertices of the once subdivided cube of side length 2. Given $E:\{0,1,2\}^k \to G$, let $F$ denote its restriction to $\{0,2\}^k$. We need to compare the currents of integration defined by the piecewise smooth interpolations $F'$ and $E'$. We construct currents $H'$ and $K'$ such that 
\begin{eqnarray*}
E'-F'=H'+\partial K'.
\end{eqnarray*}
Furthermore, $H'$ is a sum of terms, one for each face of $\{0,1,2\}^k$, which depends only on the restriction of $E$ to that face, in a functorial manner.

We use induction on $k$. When $k=1$, $E'$ has two simplices, $F'$ has one,
\begin{eqnarray*}
E'-F'&=&\sigma_{(E(0),E(1))}+\sigma_{(E(1),E(2))}-\sigma_{(F(0),F(2))}\\
&=&\sigma_{(E(0),E(1))}+\sigma_{(E(1),E(2))}-\sigma_{(E(0),E(2))}\\
&=&\partial \sigma_{(E(0),E(1),E(2))},
\end{eqnarray*}
so we set $H'=0$ and $K'=\sigma_{(E(0),E(1),E(2))}$, which can be viewed as the cone with apex $E(1)$ over $F'$.

\begin{center}
\includegraphics[width=5in]{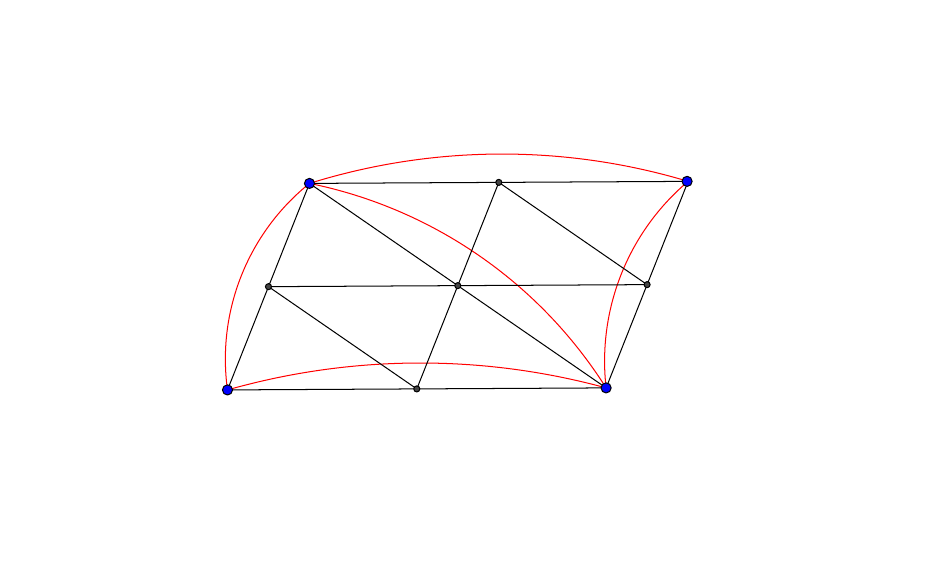}

Currents $F'$ (red) and $E'$ (black) when $k=2$. 

Current $H'$ consists of one triangle per side, one of them is thickened.
\end{center}

Let $k\geq 2$. A $k$-dimensional cube has $2k$ faces, indexed by $i\in\{1,\ldots,k\}$ and $x\in\{0,1\}$ (this describes the face whose vertices have an $x$ as $i$-th coordinate). The current $\partial E'-\partial F'$ is a signed sum of simplices which splits into $2k$ sub-sums, according to faces of $\{0,1,2\}^k$, 
\begin{eqnarray*}
\partial E'-\partial F'=\sum_{i=1}^{k}\sum_{x=0}^{1} E'_{i,x}-F'_{i,x}.
\end{eqnarray*}
Note that $E'_{i,x}$ (resp. $F'_{i,x}$) is the piecewise smooth approximation of the restriction of $E$ (resp. $F$) to the $(i,x)$-face. By induction, 
\begin{eqnarray*}
E'_{i,x}-F'_{i,x}=H'_{i,x}+\partial K'_{i,x}.
\end{eqnarray*}
By induction, $H'_{i,x}$ is a sum of terms $H'_{(i,x),(j,y)}$, one for each face of $\{0,1,2\}^{k-1}$. Each codimension 2 face $(i,x),(j,y)$ arises twice in the boundary of the $k$-cube, as a piece of the boundary of two faces, $(i,x)$ and $(j,y)$. By assumption, $H'_{(i,x),(j,y)}$ depends only on the restriction of $E$ to codimension 2 face $(i,x),(j,y)$, hence $H'_{(i,x),(j,y)}=-H'_{(i,y),(i,x)}$. Therefore, when summing up all $H'_{i,x}$, the result adds up to 0. Hence
\begin{eqnarray*}
E'_{i,x}-F'_{i,x}=\partial K'_{i,x}.
\end{eqnarray*}
If we set 
\begin{eqnarray*}
H'=\sum_{i=1}^{k}\sum_{x=0}^{1} K'_{i,x},
\end{eqnarray*}
then $\partial E'-\partial F'=\partial H'$. Note that $H'$ is indeed a sum of terms, one for each face, depending only on the restriction of $E$ to that face.

We introduce the cone $K'$ with vertex $E(1,\ldots,1)$ over $E'-F'-H'$. It is a $k+1$-dimensional current such that
\begin{eqnarray*}
E'-F'-H'=\partial K',
\end{eqnarray*}
as announced.

When passing from $F_{j}$ to $F_{j+1}$, one sums up a contribution $H'$ and $K'$ for each tiny cube, getting currents $H'_{j}$ and $K'_{j}$ which satisfy
\begin{eqnarray*}
F_{j+1}-F_j=H'_j+\partial K'_j.
\end{eqnarray*}
A cancellation takes place again: the contributions in $H'_{j}$ corresponding to a face which is shared by two tiny cubes cancel each other. Only contributions from outer faces remain. Therefore the number of terms in $H'_{j}$ is $2^{j(k-1)+1}k$, hence the number of simplices in $H'_j$ is $O(2^{j(k-1)})$. The number of simplices in $K'_j$ is $O(2^{jk})$.

If $C$ is a cycle made by piecing together $k$-cubes (for instance, the boundary of a $k+1$-cube), there are no outer faces at all, and $H'_j=0$. Note that the standard $k$-simplex has a natural subdivision into $k+1$ cubes (in barycentric coordinates, the $i$-th cube is defined by inequalities $\{x_i\geq \max_{j\not=i}x_j\}$). So this remark applies to every oriented pseudo-manifold (simplicial complex where each codimension 1 face is contained in exactly 2 faces) embedded in $\R^N$.

\section{Integrals of smooth forms along H\"older continuous maps}
\label{integral}

\subsection{Weights of differential forms in Carnot groups}

On the Lie algebra, exterior forms are graded according to \emph{weight}, i.e. eigenvalue under the 1-parameter group of dilations: let 
\begin{eqnarray*}
\Lambda^{k,w}=\{\omega\in\Lambda^{k}\mathfrak{g}^*\,;\,\delta_t^*\omega=t^{w}\omega\}, \quad \Lambda^{k,\geq w}=\bigoplus_{w'\geq w}\Lambda^{k,w'}.
\end{eqnarray*}
This filtration induces a filtration of differential forms on $G$: $\Omega^{k,\geq w}$ is the space of differential forms which pointwise belong to $\Lambda^{k,\geq w}$.

The following Lemma is merely a compactness property of straight simplices up to rescaling.

\begin{lem}\label{weight}
Let $\omega\in\Omega^{k,\geq w}$ be a smooth differential form on $G$, of weight $\geq w$. There exists a constant $C$ such that, for all sequences $(p_0,\ldots,p_k)$ of points of $G$ with diameter $\leq 1$,
\begin{eqnarray*}
|\int_{\sigma_{(p_0,\ldots,p_k)}}\omega|\leq C\,\mathrm{diameter}(\{p_0,\ldots,p_k\})^{w}\|\omega\|_{\infty}.
\end{eqnarray*}
\end{lem}

\begin{pf}
Without loss of generality, $p_0$ is the origin in $G$. Let $t=\mathrm{diameter}(\{p_0,\ldots,p_k\})$. By assumption, $t\leq 1$. Let $p'_i=\delta_{1/t}p_i$. Since
\begin{eqnarray*}
\sigma_{(p_0,\ldots,p_k)}=\delta_t \circ\sigma_{(p'_0,\ldots,p'_k)},
\end{eqnarray*}
\begin{eqnarray*}
|\int_{\sigma_{(p_0,\ldots,p_k)}}\omega|
&=&|\int_{\sigma_{(p'_0,\ldots,p'_k)}}\delta_t^*\omega|\\
&\leq& \|\sigma_{(p'_0,\ldots,p'_k)}\|_{C^1}\,t^{w}\|\omega\|_{\infty},
\end{eqnarray*}
and the derivatives of maps $\sigma_{(p'_0,\ldots,p'_k)}$ stay uniformly bounded while $\{p'_0,\ldots,p'_k\}$ remains in the unit ball of $G$.
\end{pf}

\subsection{Defining integrals of smooth forms}

\begin{thm}\label{2}
Let $\omega$ be a smooth $k$-form on $G$, of weight $\geq w$, such that $d\omega$ has weight $\geq w'$. Let $F:[0,1]^k\to G$ be $C^{\alpha}$ with $\alpha>\max\{\frac{k-1}{w},\frac{k}{w'}\}$. For every cube $Q\subset [0,1]^k$ of diameter $\delta$,
\begin{eqnarray*}
\textrm{the limit }\quad J(F(Q),\omega)=\lim_{j\to\infty}\int_{Q}F_j^{*}\omega \quad \textrm{ exists}, 
\end{eqnarray*}
and
\begin{eqnarray}\label{1}
|J(F(Q),\omega)|\leq \|F\|_{C^\alpha}^{w'} \delta^{w'\alpha}\|d\omega\|_{\infty}+\|F\|_{C^\alpha}^w \delta^{w\alpha}\|\omega\|_{\infty}.
\end{eqnarray}

Furthermore, if $\Psi$ is a closed pseudomanifold, then $J(F(\Psi),\omega)$ makes sense as soon as $\alpha>\frac{k}{w'}$. In addition, Stokes Theorem holds. If $\Omega$ is a pseudomanifold with boundary,
\begin{eqnarray*}
J(F(\Omega),d\omega)=J(F(\partial\Omega),\omega).
\end{eqnarray*}

Finally, assume that $\alpha>\max\{\frac{k}{w},\frac{k}{w'}\}$. Then
\begin{eqnarray*}
J(F(Q),\omega)=0. 
\end{eqnarray*}
\end{thm}

\begin{pf}
For each tiny cube, of diameter $O(2^{-j}\delta)$, the corresponding simplices in $G$ have diameter $O((2^{-j}\delta)^\alpha \|F\|_{C^\alpha})$.
Since $F_{j+1}-F_j=H'_j+\partial K'_j$,
\begin{eqnarray*}
|\int_{F_{j+1}(Q)}\omega-\int_{F_{j}(Q)}\omega|
&\leq& C\,(|\int_{H'_{j}}\omega|+|\int_{K'_{j}}d\omega|)\\
&\leq& C\,(\|\omega\|_\infty 2^{j(k-1)}(2^{-j}\delta)^{w\alpha}\|F\|_{C^\alpha}^w\\
&&+\|d\omega\|_\infty 2^{jk}(2^{-j}\delta)^{w'\alpha}\|F\|_{C^\alpha}^{w'}).
\end{eqnarray*}
If $\alpha>\max\{\frac{k-1}{w},\frac{k}{w'}\}$, the series converges, so the sequence $\int_{F_{j}(Q)}\omega$ has a limit $J(F(Q),\omega)$.

The direct estimate for $F_j$ is a bit weaker, since the total number of simplices is $O(2^{jk})$,  
\begin{eqnarray*}
|\int_{F_{j}(Q)}\omega|
&\leq& C\,\|\omega\|_\infty 2^{jk}(2^{-j}\delta)^{w\alpha}\|F\|_{C^\alpha}^w\\
&=& C\,\|\omega\|_\infty \|F\|_{C^\alpha}^w \delta^{w\alpha} 2^{j(k-w\alpha)}.
\end{eqnarray*}
Taking $j=0$ in the direct estimate and summing yields
\begin{eqnarray*}
|J(F(Q),\omega)|\leq C\,\max\{\|\omega\|_\infty\delta^{w\alpha}\|F\|_{C^\alpha}^w,\|d\omega\|_\infty\delta^{w'\alpha}\|F\|_{C^\alpha}^{w'}\}.
\end{eqnarray*}

For a closed pseudomanifold, $F_{j+1}(\Psi)-F_j(\Psi)=\partial K'_j$,
\begin{eqnarray*}
|\int_{F_{j+1}(\Psi)}\omega-\int_{F_{j}(\Psi)}\omega|
&\leq& C\,|\int_{K'_{j}}d\omega|\\
&\leq& C\,\|d\omega\|_\infty 2^{jk}(2^{-j}\delta)^{w'\alpha}\|F\|_{C^\alpha}^{w'},
\end{eqnarray*}
hence the series converges, and $J(F(\Psi),\omega)$ is well-defined, as soon as $\alpha>\frac{k}{w'}$. In case $\Psi=\partial\Omega$, Stokes Theorem holds for piecewise smooth interpolations, and one can pass to the limit.

The direct estimate implies that if $\alpha>\frac{k}{w}$, then $\int_{F_{j}(Q)}\omega$ tends to 0 and thus $J(F(Q),\omega)=0$.
\end{pf}

\subsection{Carnot manifolds}
\label{mani}

Since they rely merely on the filtration of differential forms by conditions weight $\geq w$, the previous considerations extend to a class of curved Carnot manifolds. By a \emph{Carnot manifold}, we mean a smooth manifold $M$ equipped with a smooth sub-bundle $\Delta\subset TM$, subject to the following non-integrability condition. Define inductively $\Delta^1=\Delta$ and $\Delta^{j+1}(p)$ as linear span the set of values at $p$ of brackets $[X,Y]$ of sections $X$ of $\Delta$ and $Y$ of $\Delta^j$. We require that there exist $s$, the \emph{step} of $\Delta$, such that $\Delta^s=TM$. 

A Carnot manifold is \emph{equiregular} if all $\Delta^j$ are smooth sub-bundles. They determine a filtration of the Lie algebra of smooth vectorfields, whence an associated graded Lie algebra $\mathcal{G}$. In the equiregular case, $\mathcal{G}$ identifies with the space of smooth sections of a smooth bundle of $n$-dimensional graded Lie algebras, the \emph{nilpotentization} $\mathfrak{n}$.

A smooth Euclidean structure on the fibers of $\Delta$ gives rise to a Carnot-Carath\'eodory metric $d$ on $M$. Any two Carnot-Carath\'eodory metrics are equivalent.

Equiregular Carnot manifolds have \emph{privileged} coordinate systems at every point $p_0$ (\cite{Bellaiche}, section 5.1). This is a local diffeomorphism $\phi:U\to \mathfrak{n}(p_0)$, $U$ a neighborhood of $p_0$ in $M$, such that if $z_j$ is a weight $w$ coordinate on $\mathfrak{g}(p_0)$, then $z_j\circ\phi=O(d(\cdot,p_0)^w)$. Say a differential form $\omega$ on $M$ has weight $\geq w$ if it is so in some privileged coordinate system. This does not depend on the choice of privileged coordinate system. Therefore the space $\Omega^{k,\geq w}$ of $k$-forms of weight $\geq w$ is well defined.

However, straight simplices do depend on the choice of coordinates, so once and for all, we make a smooth choice $(p_0,p)\mapsto (p_0,\phi_{p_0}(p))$, $V\to \mathfrak{n}$, defined on a neighborhood $V$ of the diagonal in $M\times M$. Let us define inductively
\begin{eqnarray*}
\sigma_{(p_0,\ldots,p_k)}=\phi_{p_0}^{-1}\circ\kappa_{p_0}(\phi_{p_0}\circ\sigma_{(p_1,\ldots,p_k)}),
\end{eqnarray*}
where, given a map $f$ from $k-1$-simplex to $\mathfrak{n}(p_0)$, $\kappa_{p_0}(f)$ denotes the affine cone in $\mathfrak{n}(p_0)$ with vertex at the origin over $f$.
This allows to interpolate maps from a cube to $M$. 

When a cube is subdivided, a simplex may occur in two neighboring sub-cubes. As verified in subsection \ref{interpolation}, the ordering of vertices is the same in both sub-cubes, so the parametrized simplex is the same. Therefore the cancellation mechanisms observed above persist.

According to the Ball-box Theorem, if $\phi$ is a privileged coordinate system at $p_0$, the function $p\mapsto d(p_0,\phi^{-1}(z))$ on $\mathfrak{n}(p_0)$ is equivalent to a homogeneous norm. Hence the Carnot-Carath\'eodory diameter of $(p_0,\ldots,p_k)$ and the radius of the smallest homogeneous ball containing $\phi_{p_0}^{-1}(p_1,\ldots,p_k)$ are of the same order. Hence Lemma \ref{weight} generalizes with hardly any change.

\begin{lem}\label{weightbis}
Fix a compact subset $K$ of $M$. Let $\omega\in\Omega^{k,\geq w}$ be a smooth differential form on $M$, of weight $\geq w$. There exists a constant $C$ such that, for all sequences $(p_0,\ldots,p_k)$ of points of $K$,
\begin{eqnarray*}
|\int_{\sigma_{(p_0,\ldots,p_k)}}\omega|\leq C\,\mathrm{diameter}(\{p_0,\ldots,p_k\})^{w}\|\omega\|_{\infty}.
\end{eqnarray*}
\end{lem}

\section{Upper bounds on H\"older exponents}
\label{bounds}

\subsection{Representing homology classes}

To draw conclusions from vanishing of integrals, we borrow arguments from Gromov and Rumin, see \cite{Trento}. 

Let $\mathfrak{g}$ be a Carnot Lie algebra. In each degree $k$, the Lie algebra cohomology $H^k(\mathfrak{g})$ is graded according to weights,
\begin{eqnarray*}
H^k(\mathfrak{g})=\bigoplus_{w=w(k)}^{W(k)}H^{k,w}(\mathfrak{g}).
\end{eqnarray*}

\begin{defi}
Say a Carnot manifold is \emph{equihomological} if the Hochschild cohomology spaces $\mathcal{E}_{0}(p):=H^{\cdot}(\mathfrak{n}(p))$ of its nilpotentization form a smooth bundle. 
\end{defi}

\begin{prop}
\label{top}
Let $M$ be an equiregular and equihomological Carnot manifold. Let $w(k)$ denote the minimum weight of the cohomology in degree $k$ of its nilpotentization at some point. Let $U\subset\R^n$ be a neighborhood of this point and $F:U\to V\subset G$ be a $C^\alpha$ map, with $\alpha>\frac{k}{w(k)}$. Then $F$ induces the 0 map $H^k(V,\R)\to H^k(U,\R)$.
\end{prop}

\begin{pf}
According to Rumin, \cite{Rumin}, every real cohomology class of $V$ contains a smooth differential form of weight $\geq w(k)$. Every real homology class in $U$ can be represented by an oriented pseudomanifold $X$ and a Lipschitz map $f:X\to U$, \cite{Baas}. Since $\omega$ is closed, if $\tilde{F}$ is a piecewise smooth approximation of $F$ on a neighborhood of $f(X)$, $\int_{f(X)}\tilde{F}^{*}\omega$ only depends on the homotopy class of $\tilde{F}$. Therefore this integral is equal to $J(F(f(X)),\omega)$. According to Theorem \ref{2}, this number vanishes. This shows that $\langle f(X),F^{*}[\omega]\rangle=0$. One concludes that $F^{*}[\omega]=0$, and thus that the morphism $F^*$ vanishes in degree $k$.
\end{pf}

\begin{cor}
\label{cor}
If there exists a $C^\alpha$ homeomorphism between open sets of $\R^n$ and $G$, then $\alpha\leq\min\{\frac{k}{w(k)}\,;\,k=1,\ldots,n-1\}$.
\end{cor}

\begin{pf}
Given a homeomorphism $F$ between open sets, pick a ball contained in its domain and remove an embedded $n-k-1$-sphere from the ball to get $U$. Set $V=F(U)$. By Alexander duality, $H^{k}(U,\R)\not=0$. By assumption, $F$ induces an isomorphism $H^{k}(V,\R)\to H^{k}(U,\R)$. Proposition \ref{top} implies that $\alpha\leq\frac{k}{w(k)}$.
\end{pf}

\medskip

The following cohomology calculations are borrowed from \cite{Rumin} and \cite{Trento}, section 9. Together with Corollary \ref{cor}, they provide a proof of Theorem \ref{thm}.

\begin{exa}
For all Carnot groups of Hausdorff dimension $Q$, $k=n-1$ gives $w(k)=Q-1$, whence the bound $\alpha(M)\leq\frac{n-1}{Q-1}$ for all equiregular equihomological Carnot manifolds of Hausdorff dimension $Q$.
\end{exa}

\begin{exa}
For the $4m+3$-dimensional quaternionic Heisenberg group, $k=n-2$ gives $w(k)=Q-2$, whence the bound $\alpha(G)\leq\frac{n-2}{Q-2}=\frac{4m+1}{4m+4}$.
\end{exa}

\begin{exa}
If $G$ is $2m+1$-dimensional Heisenberg group $\mathbb{H}^m$, then $k=m+1$ gives $w(k)=m+2$, whence $\alpha(M)\leq \frac{m+1}{m+2}$ for all $2m+1$-dimensional contact manifolds. 
\end{exa}

\begin{exa}
If $G$ is a Carnot group whose first layer has dimension $h$ and contains a regular isotropic $k$-plane, $w(k)=Q-k$, where Hausdorff dimension $Q=n+h$. If $h-k\geq (n-h)k$, the nilpotentization of a generic rank $h$ distribution in dimension $n$ has this property. It is also equiregular and equihomological, whence $\alpha(M)\leq \frac{n-k}{Q-k}=\frac{n-k}{n+h-k}$ for such Carnot manifolds. 
\end{exa}

%\subsection{Volumes}

%Is there a simpler argument for $\alpha\leq \frac{n-1}{Q-1}$ ? 

%Let $Q$ be a cube contained in the domain of $F$. Let $\omega$ denote a differential form of weight $Q-1$ such that $d\omega$ is the left invariant volume form on $G$. For piecewise smooth approximations,\begin{eqnarray*}\int_{F_j(Q)}d\omega=\int_{G}I(F_j(\partial Q),p)\,dp\end{eqnarray*}is the integral of the index function. Since $F_j$ converges uniformly to $F$, the index function converges away from $F(\partial Q)$. Assume that $\alpha>\frac{n-1}{Q-1}$. Then $F(\partial Q)$ has vanishing Lebesgue measure, so convergence holds almost everywhere. If the index function were dominated, one could conclude that, as $j$ tends to infinity, the right hand side converges to $\int_{G}I(F(\partial Q),p)\,dp=\mathrm{volume}(F(Q))>0$. On the other hand, if $\alpha> \frac{n-1}{Q-1}$, $\int_{F_j(Q)}d\omega=\int_{F_j(\partial Q)}\omega$ tends to $J(F(\partial Q),\omega)=0$.

%This argument seems to have a flaw: does domination hold ?

\section{Rumin currents in dimension 3}
\label{sharp}

In this section, we focus on $C^\alpha$ surfaces in the 3-dimensional Heisenberg group $G=\mathbb{H}^1$, for $\alpha>\frac{1}{2}$. Influenced by \cite{Zusttrees}, we expect that such maps factor through trees. Therefore integrals of all smooth 2-forms should vanish on them.

We first show that it suffices to prove that integrals of the two left-invariant weight 3 2-forms vanish. One can in fact express integrals of general 2-forms in terms of integrals of left-invariant forms, via Towghi integrals.

Second, we introduce a class of integral currents, called Rumin-flat integral currents, which is adapted to the anisotropic geometry of Heisenberg group. We show that the key to the sharp bound $\frac{1}{2}$ for H\"older exponents could be a slicing theorem in this class.

\subsection{Notation}

$\mathbb{H}^1$ is the set of $3\times 3$ unipotent matrices $\begin{pmatrix}
1 & x & z \\
0 &  1 & y \\
0 & 0 & 1
\end{pmatrix}$, where $x,y,z\in \R$. The forms $dx$, $dy$ and $\theta=dz-x\,dy$ constitute a basis of left-invariant 1-forms, with $dx$ and $dy$ of weight 1 and $\theta$ of weight 2. Left-invariant 2-forms are spanned by $dx\wedge dy=-d\theta$, $dx\wedge\theta$ and $dy\wedge \theta$, with $dx\wedge dy$ of weight 2, $dx\wedge\theta$ and $dy\wedge \theta$ of weight 3.

\subsection{Reduction to integrating left-invariant 2-forms}

The following merely uses the fact that $F(S)$ is a normal current satisfying
\begin{eqnarray*}
|J(F(S),\omega)|\leq C_1\|\omega\|_{\infty}\mathrm{diameter}(F(S))^{w(\omega)\alpha}+C_2\|d\omega\|_{\infty}\mathrm{diameter}(F(S))^{w(d\omega)\alpha},
\end{eqnarray*}
where the sup norm is taken on the convex hull of $F(S)$.

\begin{lem}\label{variable}
Let $S$ be the unit square, let $F:S\to\mathbb{H}^1$ be $C^\alpha$ with $\alpha>1/2$. Assume that for all sub-squares $B\subset S$ and all left-invariant weight 3 2-forms $\omega_0$ on $\mathbb{H}^1$, $J(F(B),\omega_0)=0$. Then $J(F(S),\omega)=0$ for all smooth 2-forms $\omega$.
\end{lem}

\begin{pf}
Let $\omega$ be a smooth weight 3 2-form. Subdivide $S$ in $2^{2j}$ sub-squares $B$ of diameter $2^{-j}$. For every such $B$, let $\omega_B$ be the left-invariant weight 3 2-form which coincides with $\omega$ at the image of the center of $B$. Then 
\begin{eqnarray*}
|\omega-\omega_B|\leq C\,\mathrm{diameter}(B)^{\alpha},\quad|d(\omega-\omega_B)|\leq C\,\mathrm{diameter}(B)^{\alpha}
\end{eqnarray*}
on the convex hull of $F(B)$. Since, by assumption, $J(F(B),\omega_B)=0$, 
\begin{eqnarray*}
|J(F(B),\omega)|&\leq& C_1\,\|\omega-\omega_B\|_{L^{\infty}(F(B))}\mathrm{diameter}(B)^{3\alpha}\\
&&+C_2\,\|d(\omega-\omega_B)\|_{L^{\infty}(F(B))}\mathrm{diameter}(B)^{4\alpha}\\
&\leq& C''\,\mathrm{diameter}(B)^{4\alpha-2}\mathrm{area}(B).\end{eqnarray*}
Summing up gives a bound on $J(F(S),\omega)$ that tends to 0.

Every smooth weight 2 2-form can be written $\omega=\lambda d\theta$. Since $J(F(S),d(\lambda\theta))=J(F(\partial S),\lambda\theta)=0$,
\begin{eqnarray*}
J(F(S),\lambda d\theta)=-J(F(S),d\lambda\wedge\theta)=0.
\end{eqnarray*}
Therefore integrals of all smooth 2-forms vanish.
\end{pf}

\begin{rem}
\label{3=>2}
Roger Z\"ust explained to us that if we can prove that $J(F(S),\omega)=0$ for all smooth closed weight 2 2-forms $\omega$, then $F$ cannot be injective on $S$ (in fact, $F$ factors through a tree, see \cite{Zusttrees}).
\end{rem}

\subsection{Relation to Towghi's integrals}

\begin{lem}\label{towghi}
Let $\omega$ be a smooth weight 3 2-form and $f$ a smooth function on $\mathbb{H}^1$. 
Let $\omega'=f\omega$. 
Let $F:S\to\mathbb{H}^1$ be a $C^\alpha$ map from a square to $\mathbb{H}^1$, $\alpha>\frac{1}{2}$. 
Then $J(F(S),\omega')$ can be expressed as a Towghi integral as follows, \cite{Towghi}.
For $(s,t)\in S$, let $S(s,t)\subset S$ denote the rectangle with vertices $(0,0)$ and $(s,t)$. Define
\begin{eqnarray*}
g(s,t)=J(F(S(s,t)),\omega).
\end{eqnarray*} 
Then
\begin{eqnarray*}
J(F(S),\omega')=\int_{S}(f\circ F)\,dg
\end{eqnarray*}
in Towghi's sense.
\end{lem}

\begin{pf}
Recall the definition of Towghi's integral $\int_{S}f\,dg$. Pick a subdivision $\pi$ of $S$, i.e. a pair of subdivisions $\{s_i\}$ and $\{t_j\}$ of the sides of $S$. 
Pick points $(\eta_i,\nu_j)\in[s_{i-1},s_i]\times[t_{j-1},t_j]$. Consider sums
\begin{eqnarray*}
L(f,g,\pi)=\sum_{i,j}f(\eta_i,\nu_j)\left(g(s_i,t_j)-g(s_{i-1},t_j)-g(s_i,t_{j-1})+g(s_{i-1},t_{j-1})\right).
\end{eqnarray*}
Under suitable assumptions on $f$ and $g$, these sums converge to a number which is $\int_{S}f\,dg$ by definition.

Subdivide $S$ into equal squares $B$ and pick a point $p_B$ in each of them. As in Lemma \ref{variable}, 
\begin{eqnarray*}
|\omega'-f(F(p_B))\omega|\leq C\,\mathrm{diameter}(B)^{\alpha},\quad|d(\omega'-f(F(p_B))\omega)|\leq C\,\mathrm{diameter}(B)^{\alpha}
\end{eqnarray*}
on the convex hull of $F(B)$. Therefore
\begin{eqnarray*}
|J(F(B),\omega')&-&f(F(p_B))J(F(B),\omega')|\\
&\leq& C_1\,\|\omega'-f(F(p_B))\omega\|_{L^{\infty}(f(B))}\mathrm{diameter}(B)^{3\alpha}\\
&&+C_2\,\|d(\omega'-f(F(p_B))\omega)\|_{L^{\infty}(f(B))}\mathrm{diameter}(B)^{4\alpha}\\
&\leq& C''\,\mathrm{diameter}(B)^{4\alpha-2}\mathrm{area}(B).\end{eqnarray*}
Summing up gives
\begin{eqnarray*}
J(F(S),\omega')=\sum_{B\subset S}f(F(p_B))J(F(B),\omega)+O(\mathrm{diameter}(B)^{4\alpha-2})\mathrm{area}(S).
\end{eqnarray*}
Write $B=[s_{i-1},s_i]\times[t_{j-1},t_j]$. Then
\begin{eqnarray*}
J(F(B),\omega)=g(s_i,t_j)-g(s_{i-1},t_j)-g(s_i,t_{j-1})-g(s_{i-1},t_{j-1}).
\end{eqnarray*}
So if Towghi's sums converge, their limit is equal to $J(F(S),\omega')$.

Towghi requires assumptions on $(p,q)$-variations. The $(p,p)$-variation of $g$ is 
\begin{eqnarray*}
\sup_{\pi}\sum_{i,j}|g(s_i,t_j)-g(s_{i-1},t_j)-g(s_i,t_{j-1})-g(s_{i-1},t_{j-1})|^p
\end{eqnarray*}
In the special case of a dyadic subdivision, we know from Lemma \ref{2} that 
\begin{eqnarray*}
|J(F(B),\omega)|=O(\mathrm{diameter}(B)^{3\alpha}).
\end{eqnarray*}
This gives
\begin{eqnarray*}
\sum_{B\subset F}|J(F(B),\omega)|^{p}=O(\mathrm{diameter}(B)^{3p\alpha-2})\mathrm{area}(S).
\end{eqnarray*}
So $g$ has finite (dyadic) $(p,p)$-variation for $p=\frac{2}{3\alpha}$. 
For $f\circ F$, we know that the contribution of a small square $B$ is at most $\mathrm{diameter}(B)^{q\alpha}$, so the sum of $q$-powers is at most $\mathrm{diameter}(B)^{q\alpha-2}\mathrm{area}(S)$. 
Thus $f\circ F$ has finite (dyadic) $(q,q)$-variation for $q=\frac{2}{\alpha}$. Since
\begin{eqnarray*}
\frac{1}{p}+\frac{1}{q}=\frac{3\alpha}{2}+\frac{\alpha}{2}=2\alpha>1,
\end{eqnarray*}
Towghi's assumption is satisfied for dyadic subdivisions, so Towghi's sums indeed converge.
\end{pf}

\begin{rem}
It not so clear wether the $(p,p)$ and $(q,q)$ variations are finite in unrestricted sense. 
\end{rem}

Here is a weight 2 version of Lemma \ref{towghi}.

\begin{lem}
Let $\omega$ be a smooth weight 2 2-form. Let $F:S\to\mathbb{H}^1$ be a $C^\alpha$ map from a square to $\mathbb{H}^1$, $\alpha>\frac{1}{2}$. 
Then $J(F(S),\omega)$ can be expressed as a sum of two Towghi integrals as follows.
For $(s,t)\in S$, let $S(s,t)\subset S$ denote the rectangle with vertices $(0,0)$ and $(s,t)$. Define
\begin{eqnarray*}
h(s,t)=J(F(S(s,t)),x\,dx\wedge dy),\quad k(s,t)=J(F(S(s,t)),y\,dx\wedge dy).
\end{eqnarray*} 
Write $\omega=a\,dx\wedge dy$. Then
\begin{eqnarray*}
J(F(S),\omega)=\int_{S}^T (Xa\circ F)\,dh+\int_{S}^T (Ya\circ F)\,dk
\end{eqnarray*}
in Towghi's sense.
\end{lem}

\begin{pf}
Subdivide $S$ into squares $B$. For each square $B\subset S$ with center $p_B$, let 
\begin{eqnarray*}
a&=&a(p_B)+(Xa)(F(p_B))(x-x_{F(p_B)})+(Ya)(F(p_B))(y-y_{F(p_B)})\\&&+O((x-x_{F(p_B)})^2+(y-y_{F(p_B)})^2+|(z-z_{F(p_B)}|)
\end{eqnarray*} 
be the weight 1 Taylor expansion of $a$ at $F(p_B)$. Set 
\begin{eqnarray*}
\omega_B=(a(p_B)+(Xa)(F(p_B))(x-x_{F(p_B)})+(Ya)(F(p_B))(y-y_{F(p_B)}))dx\wedge dy. 
\end{eqnarray*}
Then, on the convex hull of $F(B)$,
\begin{eqnarray*}
|\omega-\omega_B|\leq \mathrm{diameter}(B)^{2\alpha},\quad |d(\omega-\omega_B)|\leq \mathrm{diameter}(B)^{2\alpha}.
\end{eqnarray*}
Since, according to Lemma \ref{2},
\begin{eqnarray*}
|J(F(B),\omega)-J(F(B),\omega_B)|&\leq& C_1\,\|\omega-\omega_B\|_{L^{\infty}(f(B))}\mathrm{diameter}(B)^{2\alpha}\\
&&+C_2\,\|d(\omega-\omega_B)\|_{L^{\infty}(F(B))}\mathrm{diameter}(B)^{4\alpha}\\
&\leq& C''\,\mathrm{diameter}(B)^{4\alpha-2}\mathrm{area}(B),\end{eqnarray*}
\begin{eqnarray*}
J(F(S),\omega)&=&\sum_{B\subset S}J(F(B),\omega)\\
&=&\sum_{B\subset S}J(F(B),\omega_B)+O(\mathrm{diameter}(B)^{4\alpha-2})\mathrm{area}(S).
\end{eqnarray*}
Since $J(F(B),dx\wedge dy)=-J(F(B),d\theta)=0$ for all $B$,
\begin{eqnarray*}
J(F(B),\omega_B)&=&a(p_B)J(F(B),dx\wedge dy)\\
&&+(Xa)(F(p_B))J(F(B),(x-x_{F(p_B)})dx\wedge dy)\\
&&+(Ya)(F(p_B))J(F(B),(y-y_{F(p_B)})dx\wedge dy)\\
&=&(Xa)(F(p_B))J(F(B),x\,dx\wedge dy)\\
&&+(Ya)(F(p_B))J(F(B),y\,dx\wedge dy).
\end{eqnarray*}
This yields the expression of $J$ as a sum of Towghi integrals. These integrals make sense since 
\begin{eqnarray*}
J(F(B),x\,dx\wedge dy)=J(F(B),(x-x_{F(p_B)})dx\wedge dy)=O(\mathrm{diameter}(B)^{3\alpha}),
\end{eqnarray*}
which shows that $h$ has finite (dyadic) $(p,p)$-variation for $p=\frac{2}{3\alpha}$.

Alternative (shorter) argument, based on Lemma \ref{towghi}. From $d(-a\theta)=a\,dx\wedge dy -da\wedge\theta$ and $J(F(S),a\theta)=0$, get
\begin{eqnarray*}
J(F(S),a\,dx\wedge dy)&=&J(F(S),da\wedge\theta)\\
&=&J(F(S),Xa\,dx\wedge\theta)+J(F(S),Ya\,dy\wedge\theta)\\
&=&\int_{S}(Xa\circ F)\,dh'+\int_{S}(Ya\circ F)\,dk',
\end{eqnarray*}
where $h'(s,t)=J(F(S(s,t)),dx\wedge\theta)=J(F(S(s,t)),x\,dx\wedge dy)=h(s,t)$.
\end{pf}

%\subsection{Relating Towghi integrals to Z\"ust-like integrals}

%\begin{rem}The map $\Phi\mapsto J(F(Q),\Phi(x,y)\,dx\wedge dy)$ is a Radon measure on $C^0(\R^2)$.\end{rem}Indeed, since $d(\Phi(x,y)\,dx\wedge dy)=0$, the direct estimate shows that $J(F(Q),$ $\Phi(x,y)\,dx\wedge dy)$ is controlled by $\|\Phi\|_\infty$ only. This is in contrast with Towghi's expression in terms of derivatives.

%Let $S$ be a square tiled with subsquares $B$. Pick a point $p_B\in B$. Let $F(p_B)=(x_B,y_B,z_B)$. Then, since $dx\wedge dy$ has weight 2,\begin{eqnarray*}&&|(Xa)(F(p_B))J(F(B),(x-x_B)\,dx\wedge dy)-J(F(B),(x-x_B)\,(Xa)\,dx\wedge dy)|\\&\leq& C\,\mathrm{diameter}(B)^{2\alpha}\max_{F(B)} |(x-x_B)Xa-(x-x_B)(Xa)(F(p_B))|\\&\leq& C'\,\mathrm{diameter}(B)^{4\alpha}.\end{eqnarray*}Summing up yields\begin{eqnarray*}J(F(S),\omega)&=&\sum_B J(F(B),\omega)\\&=&\sum_B  (J(F(B),((x-x_B)(Xa)+(y-y_B)(Ya))\,dx\wedge dy)+O(\mathrm{diameter}(B)^{4\alpha})).\end{eqnarray*}The error term disappears since\begin{eqnarray*}\sum_B\mathrm{diameter}(B)^{4\alpha}\leq \mathrm{diameter}(B)^{4\alpha-2}\sum_B \mathrm{area}(B)\end{eqnarray*}tends to 0.

\subsection{Horizontal straight simplices}

\begin{lem}
\label{horizontal}
Let $\Delta^k$ denote the standard $k$-simplex. There exist continuous, $\delta_t$-equivariant and translation equivariant maps $\sigma^{h}:(\mathbb{H}^1)^k\to PL(\Delta^2,\mathbb{H}^1)$, $k=0,1,2,3$, such that 
\begin{enumerate}
  \item $\sigma^h_{p_0,\ldots,p_k}$ is a singular simplex with vertices $p_0,\ldots,p_k$.
  \item the faces of $\sigma^h_{p_0,\ldots,p_k}$ are $\sigma^h_{p_0,\ldots,\widehat{p_i},\ldots,p_k}$'s.
  \item $\sigma^h_{p_0,p_1}$ is a horizontal piecewise linear curve.
\end{enumerate}
\end{lem}

\begin{pf}
Let $U\subset \mathbb{H}^1$ be a polyhedron containing the origin in its interior and intersecting each nontrivial orbit of the dilation group exactly once. Thanks to the h-principle, there is a $PL$ map 
\begin{eqnarray*}
p_1\mapsto \sigma^h_{e,p_1}
\end{eqnarray*}
from the $\partial U$ to the space of piecewise linear horizontal arcs, such that $\sigma^h_{e,p_1}$ joins $e$ to $p_1$. Dilations and left-translations extend this into a continuous family $(p_0,p_1)\mapsto\sigma^h_{p_0,p_1}$.

For $p_1$ in the $\partial U$ and $p_2$ in $U$, or if $p_2$ belongs to $\partial U$ and $p_1$ to $U$, consider the map
\begin{eqnarray*}
\partial\Delta^2\to \mathbb{H}^1
\end{eqnarray*}
defined by maps $\sigma^h_{e,p_1}$, $\sigma^h_{p_1,p_2}$ and $\sigma^h_{p_2,e}$ along sides. Extend it into a $PL$ map
\begin{eqnarray*}
\partial\Delta^2\to \mathbb{H}^1,
\end{eqnarray*}
by constructing the cone with vertex at the center of gravity $\frac{1}{3}(\log p_1+\log p_2)$ in exponential coordinates. Again, dilating and left-translating defines the continuous family $(p_0,p_1,p_2)\mapsto\sigma^h_{p_0,p_1,p_2}$.

A similar cone construction defines $(p_0,p_1,p_2,p_3)\mapsto\sigma^h_{p_0,p_1,p_2,p_3}$.
\end{pf}

It follows from Lemma \ref{horizontal} that the currents $F_j(I)$ approximating a $C^\alpha$ arc $F:[0,1]\to \mathbb{H}^1$, $\alpha>\frac{1}{2}$, can be chosen to be piecewise $PL$ horizontal arcs.

\subsection{Rumin-flat integral currents}

Rumin's complex has one weight, either $0,1,3$ or $4$ in each degree $0,1,2,3$. Rumin's differential on 0-forms is the weight 1 component of the exterior differential. Rumin's differential on weight 3 2-forms is equal to the exterior differential. The last Rumin differential maps weight 1 1-forms to weight 3 2-forms as follows. If $\omega$ is a smooth weight 1 1-form, there is a unique smooth function $f$ such that $d(\omega+f\theta)$ has weight 3, denote it by $d_R\omega$. $d_R$ is a second order differential operator. For other degrees, we use the notation $d_R$ as well.

A \emph{Rumin current} is a continuous linear functional on the space of smooth compactly supported Rumin forms. The \emph{Rumin boundary} of a Rumin current is defined by duality: for a smooth Rumin form $\phi$,
\begin{eqnarray*}
(\partial_{R}T)(\phi)=T(d_R \phi).
\end{eqnarray*}
The \emph{Rumin mass} of a Rumin current is its norm as a functional on $C^0$ Rumin forms.

A current $T$ defines a Rumin current. By convention, let us set $(\partial_{R}T)(\phi)=0$ if $\phi$ is a 1-form or a 2-form of weight 2. This completes the definition of $\partial_{R}T$ as a current. Then, on smooth 1-forms $\phi$,
\begin{eqnarray*}
(\partial_R T)(\phi)=T(d(\phi-\frac{d\phi}{d\theta}\theta))=\partial T(\phi-\frac{d\phi}{d\theta}\theta).
\end{eqnarray*}
Note that when $T$ is an integral 2-current, $\partial_R T$ need not be an integral current. It is if $\partial T$ vanishes on weight 2 1-forms, since then $\partial_R T=\partial T$. 

\begin{defi}
For a Rumin current $T$, define the \emph{Rumin-flat norm}, $\|T\|_{R\flat}$, as the best constant $C$ in the inequality
\begin{eqnarray*}
|T(\omega)|\leq C\max\{\|\omega\|_{\infty},\|d_R\omega\|_{\infty}\}.
\end{eqnarray*}
Equivalently,
\begin{eqnarray*}
\|T\|_{R\flat}=\inf\{M(R)+M(S)\,;\,R,\,S \textrm{ Rumin currents, }T=R+\partial_R S\}.
\end{eqnarray*}
\end{defi}

Since there are less test forms in the Rumin complex than in the de Rham complex, Rumin masses are less than usual masses, and Rumin-flat convergence implies flat convergence.

\begin{defi}
A piecewise $C^1$ chain in $\mathbb{H}^1$ is \emph{horizontal} if it is an integral linear combination of piecewise $C^1$ simplices whose edges (1-dimensional faces) are horizontal PL curves.

A current $T$ is a \emph{Rumin-flat integral current} if it is a Rumin-flat limit of horizontal piecewise $C^1$ chains whose simplices have diameters tending uniformly to 0. 
\end{defi}

Rumin flat integral currents are special cases of integral currents. By definition, for every Rumin-flat integral 1-current $T$, and every smooth function $\lambda$ on $\mathbb{H}^1$, $T(\lambda\theta)=0$. Note that horizontal piecewise $C^1$ 2-chains are usually not Rumin-flat integral currents. Rumin-flat integral 2-currents are really weird things. Nevertheless, for Rumin-flat integral currents, $\partial_R T=\partial T$ is integral.

\begin{prop}
$C^\alpha$ arcs or squares in $\mathbb{H}^1$, $\alpha>\frac{1}{2}$, are Rumin-flat integral currents.

Furthermore, the following variant of the direct estimate (\ref{1}) of Theorem \ref{2} holds. For a $C^\alpha$ arc $F(I)$, $I$ an interval of length $\delta$, for every weight 1 1-form $\omega$,
\begin{eqnarray}
|J(F(I),\omega)|\leq \|F\|_{C^\alpha} \delta^{\alpha}\|\omega\|_{\infty}+\|F\|_{C^\alpha}^{3} \delta^{3\alpha}\|d_R\omega\|_{\infty}.
\end{eqnarray}
If $S=[0,1]^2$, for a $C^\alpha$ square $F(S)$, 
$$
\partial_R F(S)=F(\partial S).
$$
\end{prop}

\begin{pf}
Let $Q=[0,1]$ or $[0,1]^2$. Let $F:Q\to \mathbb{H}^1$ be $C^\alpha$ with $\alpha>\frac{1}{2}$. Using the horizontal simplices of Lemma \ref{horizontal} as a replacement for those of subsection \ref{straight}, we get a sequence of piecewise smooth maps $F'_j$ converging uniformly to $F$, such that
\begin{eqnarray*}
F'_{j+1}-F'_j=H''_j+\partial K''_j,
\end{eqnarray*}
where $H''_j$ and $K''_j$ are sums of horizontal simplices of small diameter. If $Q=[0,1]$, since the edges of each simplex $\sigma$ of $K''_j$ are horizontal, $\partial\sigma=\partial_R \sigma$, hence $\partial K''_j=\partial_R K''_j$. Thus, in both cases,
\begin{eqnarray*}
F'_{j+1}-F'_j=H''_j+\partial_R K''_j,
\end{eqnarray*}
Since Rumin masses are smaller than usual masses, the sequence $F'_j$ converges in Rumin-flat norm, hence $F(Q)$ is a Rumin-flat integral current. 

If $Q=S=[0,1]^2$, for every weight 1 1-form $\omega$,
\begin{eqnarray*}
J(F(\partial S),\omega)=J(F(\partial S),\omega+f\theta)=J(F(S),d(\omega+f\theta))=J(F(S),d_R\omega),
\end{eqnarray*}
showing that $\partial_R F(S)=F(\partial S)$.
\end{pf}

\subsection{Slicing and Rumin-flatness}

Federer's Slicing Theorem, \cite{Federer}, Theorem 4.3.2 page 438, states that, given a Lipschitz function $\R^{m}\to\R^n$, a flat current $T$ of dimension $k$ on $\R^{m}$ can be expressed as the integral over $\R^n$ of a family of flat currents $\langle T,f,z\rangle$, $z\in\R^n$, of dimension $k-n$, supported on fibers $f^{-1}(z)$ of $f$. One can think of $\langle T,f,z\rangle$ as the intersection of $T$ and $f^{-1}(z)$. 

\medskip

\textbf{Question}. Let $T$ be a Rumin-flat integral 2-current in $\mathbb{H}^1$ and let $f:\mathbb{H}^1\to\R$ be a Lipschitz function. Is it true that for almost every $t\in\R$, $\langle T,f,t\rangle$ is a Rumin-flat integral 1-current ? 

\medskip

We would need this merely for $C^\alpha$ squares $T=F(S)$ and $f=x$ coordinate. Since $f$ is Euclidean Lipschitz, Federer's Slicing Theorem asserts that for all real valued Baire functions $\Phi$ on $\R$, and all smooth compactly supported $1$-forms $\psi$ on $\mathbb{H}^1$,
\begin{eqnarray*}
\int_{\R}\langle F(S),f,x\rangle(\psi)\,\Phi(x)\,dx=J(F(S),\psi\wedge \Phi(x)\,dx).
\end{eqnarray*}
Take $\Phi=1$ and $\psi=\chi\, \theta$, where $\chi$ is a cut-off which equals 1 on a neighborhood of $F(S)$. This yields
\begin{eqnarray*}
\int_{\R}\langle F(S),f,x\rangle(\theta)\,dx=J(F(S),\theta\wedge dx).
\end{eqnarray*}
If almost all slices $\langle F(S),f,x\rangle$ are Rumin flat integral 1-currents, then \break $\langle F(S),f,x\rangle(\theta)=0$ for a.e. $x$, and $J(F(S),\theta\wedge dx)=0$. Applying a rotation would imply that $J(F(S),\theta\wedge dy)=0$ as well. According to Lemma \ref{variable}, this suffices in order that all integrals of smooth 2-forms vanish, i.e. the currents $F(S)$ vanish identically. This implies that $C^\alpha$ maps $\R^3\to\mathbb{H}^1$ induce 0 on the homology of arbitrary open subsets, a fact which cannot happen for homeomorphisms.

\bigskip
\noindent
Zoltan Balogh
\par\noindent Mathematisches Institut
\par\noindent Universit\"at Bern
\par\noindent Sidlerstrasse 5
\par\noindent 3012 Bern, Switzerland
\par\noindent 
e-mail: zoltan.balogh@math.unibe.ch

\bigskip

\noindent
Artem Kozhevnikov
\par\noindent tinyclues
\par\noindent https://www.tinyclues.com
\par\noindent
e-mail: artem@tinyclues.com
\bigskip

\noindent
Pierre Pansu 
\par\noindent Laboratoire de Math\'ematiques d'Orsay
\par\noindent Univ. Paris-Sud, CNRS, Universit\'e
Paris-Saclay
\par\noindent 91405 Orsay, France
\par\noindent
e-mail: pierre.pansu@math.u-psud.fr

\end{document}